\documentclass{article}
\usepackage{graphicx}
\usepackage{latexsym,amsmath}
\usepackage{amsmath,amsthm}
\usepackage{amsfonts}
\usepackage[psamsfonts]{amssymb}
\usepackage{enumerate}

\usepackage{url}

\usepackage{pstricks}
\usepackage{epstopdf}

\theoremstyle{plain}
\newtheorem{theorem}{Theorem}

\newtheorem*{theorem*}{Theorem}
\newtheorem*{dp*}{Dichotomy Principle (DP)}
\newtheorem*{cp*}{Containment Principle (CP)}
\newtheorem*{hfcp*}{Finite Containment Principle (FCP)}
\newtheorem*{hjcp*}{Jordan Filling Principle (JFP)}
\newtheorem*{DP*}{Darboux-Picard Theorem (DPT)}

\theoremstyle{definition}

\theoremstyle{definition}

\theoremstyle{remark}

\theoremstyle{remark}

\newtheorem*{remark*}{Remark}

\newcommand{\R}{\mathbb{R}}

\newcommand{\dd}{{\mathrm{\,d}}}

\renewcommand{\P}{\operatorname{\mathsf{P}}}
\newcommand{\E}{\operatorname{\mathsf{E}}}
\newcommand{\Var}{\operatorname{\mathsf{Var}}}

\renewcommand{\le}{\leqslant}
\renewcommand{\ge}{\geqslant}

\begin{document}
\thispagestyle{empty}


\title{Between Chebyshev and Cantelli}
\author{Iosif Pinelis\thanks{Supported by NSF grant DMS-0805946}}
\maketitle

\begin{abstract}
A family of exact upper bounds interpolating between Chebyshev's and Cantelli's is presented. 
\end{abstract}



Let $X$ be any zero-mean unit-variance random variable (r.v.): $\E X=0$ and $\Var X=1$. 
\big(Obviously, for any non-degenerate r.v.\ $Y$ with a finite second moment, its standardization $\frac{{Y-\E Y}}{\sqrt{\Var Y}}$ is a zero-mean unit-variance r.v.\big) 
Take any $b\in(0,\infty)$. 
Chebyshev's inequality states that 
\begin{equation}\label{eq:cheb}
	\P(|X|\ge b)\le\frac1{b^2}. 
\end{equation}
Cantelli's bound on the probabilities of one-sided deviations of $X$ from $0$ is obviously smaller: 
\begin{equation}\label{eq:cant}
	\P(X\ge b)\le\frac1{1+b^2}. 
\end{equation}
Moreover, Cantelli's bound is exact, as it is attained when $X$ takes on values $-1/b$ and $b$ with probabilities $\frac{b^2}{1+b^2}$ and $\frac1{1+b^2}$, respectively. 
Chebyshev's bound is also exact, but only for $b\ge1$: indeed, let $X$ take on values $-b$, $0$, and $b$ with probabilities $\frac1{2b^2}$, $1-\frac1{b^2}$, and $\frac1{2b^2}$, respectively. 
The obviously modified Chebyshev's bound given by the inequality 
\begin{equation}\label{eq:cheb1}
	\P(|X|\ge b)\le1\wedge\frac1{b^2}  
\end{equation}
is exact for all $b\in(0,\infty)$; indeed, for $b\in(0,1)$ let $X$ take on each of the values $\pm1$ with probability $\frac12$. 
Clearly, Cantelli's bound is still smaller than modified Chebyshev's, for all $b\in(0,\infty)$. 

Observe that the event 
$\{|X|\ge b\}$ under the probability sign in \eqref{eq:cheb1} means that $X$ takes on a value outside the symmetric interval $(-b,b)$, whereas the event 
$\{X\ge b\}$ under the probability sign in \eqref{eq:cant} means that $X$ takes on a value outside the utterly asymmetric interval $(-\infty,b)$. 
More generally, one may ask about the exact upper bound on $\P\big(X\notin(-a,b)\big)$, for any given interval $(-a,b)$ containing $0$. 
The need for such a bound, which would in this sense interpolate between Chebyshev's and Cantelli's, arises naturally in studies of the distributions of the so-called self-normalized sums \cite{BE-self-norm}, where one needs a good upper bound on that probability that a quadratic polynomial $X^2+AX+B$ in a r.v.\ $X$ will take on a nonnegative value. 

When considering the probability $\P\big(X\notin(-a,b)\big)$, without loss of generality one may assume that $a\ge b$. Indeed, $\P\big(X\notin(-a,b)\big)=\P\big(-X\notin(-b,a)\big)$, and the r.v.\ $-X$ is zero-mean and  unit-variance whenever $X$ is so. 
Accordingly, let us present 

\begin{theorem}\label{th:}
Take any $a$ and $b$ such that $0<b\le a<\infty$. 
Then
\begin{equation}\label{eq:bound}
	\P\big(X\notin(-a,b)\big)\le P_{a,b}:=
	\left\{
	\begin{alignedat}{2}
&	1 &\text{\quad if \quad}& ab\le1, \\
&	\frac{4+(a-b)^2}{(a+b)^2} &\text{\quad if \quad}& \frac{(a-b)b}2\le1\le ab, \\
&	\frac1{1+b^2} &\text{\quad if \quad}& 1\le\frac{(a-b)b}2, 	
	\end{alignedat}
	\right.
\end{equation}
and this upper bound is exact. 
Moreover, 
\begin{equation}\label{eq:always}
	\P\big(X\notin(-a,b)\big)\le 
	1\wedge\frac{4+(a-b)^2}{(a+b)^2}
\end{equation}
in all of the three cases in \eqref{eq:bound}. 
\end{theorem}

Note that $P_{b,b}$ coincides with the modified Chebyshev bound \eqref{eq:cheb1}, whereas $P_{\infty,b}:=\lim_{a\to\infty}P_{a,b}$ coincides with the Cantelli bound \eqref{eq:cant}. 
So, letting $a=kb$ and varying $k$ from $1$ to $\infty$, one obtains a decreasing family $(P_{kb,b}\colon1\le k\le\infty)$ of exact upper bounds interpolating between modified Chebyshev's and Cantelli's. 
The members of this family of bounds with $k\in\{1,\dots,6,\infty\}$ are shown in the picture here. 
One can see that even for such moderate values of the ``asymmetry parameter'' $k$ as $2$ or $3$, the improvement of the bound $P_{kb,b}$ over Chebyshev's may be quite significant; for instance, Chebyshev's bound $P_{1,1}=1$ is $80\%$ greater than $P_{2,1}=\frac59$, and it is $100\%$ greater than $P_{k,1}=\frac12$ for $k\ge3$. 

\begin{center}
\includegraphics[scale=.9]{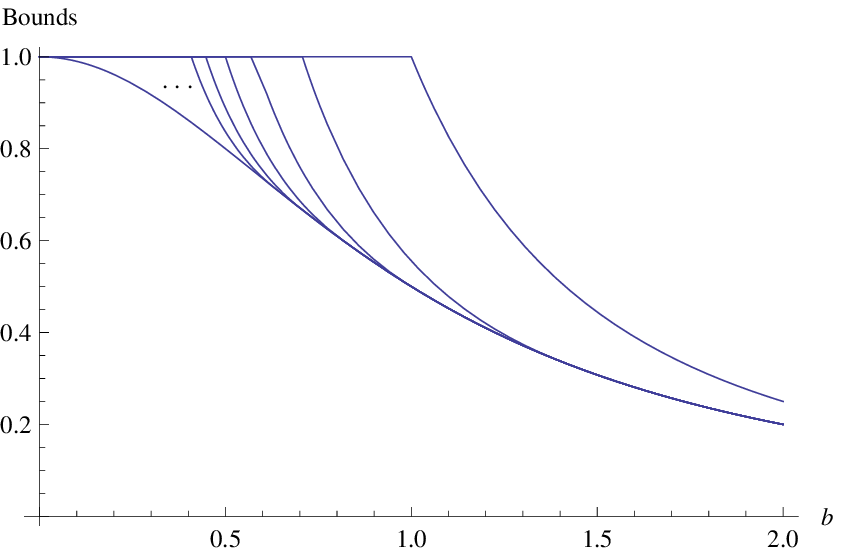}
\end{center}

Theorem~\ref{th:} can be proved by a method going back to Chebyshev and Markov; cf.\ e.g.\  \cite{krein-nudelman,karlin-studden,kemper-dual,pin98}. Rewrite $\P\big(X\notin(-a,b)\big)$ as $\int g(X)\dd\P$, where $g:=\chi_{\R\setminus(a,b)}$. Then one can try to find the best possible upper bound on $\int g(X)\dd\P$ as $\inf\int f(X)\dd\P$, where the infimum is taken over all functions $f$ that majorize $g$ 
and are linear combinations 
of the moment functions $x\mapsto1$, $x\mapsto x$, and $x\mapsto x^2$, corresponding to the restrictions $\int\dd\P=1$, $\int X\dd\P=0$, and $\int X^2\dd\P=1$. 
So, in our optimization problem the function $f$ is a quadratic polynomial such that $f\ge g$ on $\R$. 
Take now the majorizing $f(x)$ to be $\equiv1$, 
$\equiv\big(\frac{2x+a-b}{a+b}\big)^2$, or $\equiv\big(\frac{bx+1}{b^2+1}\big)^2$ 
in the three respective cases in \eqref{eq:bound};  
actually, in all of the three cases one has $g(x)\le1\wedge\big(\frac{2x+a-b}{a+b}\big)^2$ for all $x\in\R$. 
Next, writing $\P\big(X\notin(-a,b)\big)=\int g(X)\dd\P\le\int f(X)\dd\P$, and then taking into account 
the restrictions $\int\dd\P=1$, $\int X\dd\P=0$, and $\int X^2\dd\P=1$, one obtains the inequalities in \eqref{eq:bound} and \eqref{eq:always}. 
The exactness of the bound $P_{a,b}$ 
follows since it is attained when a r.v.\ $X$ takes on the values 
\begin{enumerate}[(i)]
	\item $-\sqrt{\frac ab}$  and $\sqrt{\frac ba}$ with respective probabilities $\frac b{a+b}$  and $\frac a{a+b}$  --- when $ab\le1$;  
	\item $-a$, $\frac{b-a}2$, and $b$ with respective probabilities $\frac{2-(a-b)b}{(a+b)^2}$, $\frac{4(ab-1)}{(a+b)^2}$, and $\frac{2+(a-b)a}{(a+b)^2}$ --- when $\frac{(a-b)b}2\le1\le ab$;  
	\item $-\frac1b$ and $b$ with respective probabilities $\frac{b^2}{1+b^2}$  and $\frac1{1+b^2}$ --- when $1\le\frac{(a-b)b}2$.  
\end{enumerate}

\bibliographystyle{acm}
\bibliography{C:/Users/Iosif/Documents/mtu_home01-30-10/bib_files/citations}

{\parskip0pt \parindent0pt \it Department of Mathematical Sciences

Michigan Technological University

Houghton, MI 49931

ipinelis@mtu.edu}
\end{document}